\documentclass{amsart}

\usepackage{fancyhdr}
\usepackage{lastpage}
\usepackage{stmaryrd,yhmath}
\usepackage{tikz}
\usepackage[all]{xy}

\newcommand{\bdis}{\begin{displaymath}}
\newcommand{\edis}{\end{displaymath}}
\newcommand{\be}{\begin{equation}}
\newcommand{\ee}{\end{equation}}
\newcommand{\mbb}{\mathbb}
\newcommand{\mcal}{\mathcal}

\theoremstyle{definition}

\theoremstyle{remark}
\newtheorem{remark}[]{Remark}

\newtheorem*{mydef1}{{\bf Theorem}}

\newtheorem*{mydef11}{{\bf Theorem 1}}

\newtheorem*{mydef12}{{\bf Theorem 2}}

\newtheorem*{mydef13}{{\bf Theorem 3}}

\newtheorem*{mydef2}{{\bf Definition}}

\newtheorem*{mydef52}{{\bf Lemma 2}}

\newtheorem*{mydef53}{{\bf Lemma 3}}

\newtheorem*{mydef54}{{\bf Lemma 4}}

\newtheorem*{mydefCHF1}{{\bf Complete Hybrid Formula 1}} 

\newtheorem*{mydefCHF2}{{\bf Complete Hybrid Formula 2}}  

\newtheorem*{mydefCHF3}{{\bf Complete Hybrid Formula 3}}

\numberwithin{equation}{section}

\begin{document}

\title{Jacob's ladders, crossbreeding in the set of $\zeta$-factorization formulas and selection of families of $\zeta$-kindred real continuous functions}

\author{Jan Moser}

\address{Department of Mathematical Analysis and Numerical Mathematics, Comenius University, Mlynska Dolina M105, 842 48 Bratislava, SLOVAKIA}

\email{jan.mozer@fmph.uniba.sk}


\begin{abstract}
In this paper we introduce the notion of $\zeta$-crossbreeding in a set of $\zeta$-factorization formulas and also the notion of complete hybrid formula as the final result of
that crossbreeding. The last formula is used as a criterion for selection of families of $\zeta$-kindred elements in class of real continuous functions.

Dedicated to recalling of Gregory Mendel's pea-crossbreeding.
\end{abstract}
\maketitle

\section{Introduction}

\subsection{}

Let us remind the following notions we have introduced (see \cite{1} -- \cite{7}) within the theory of the Riemann-zeta function:
\begin{itemize}
\item[(A)] Jacob's ladders, \cite{1}, comp. \cite{2},\cite{3},
\item[(B)] $\zeta$-oscillating system, \cite{7}, (1.1),
\item[(C)] $\zeta$-factorization formula, \cite{5}, (4.3) -- (4.18), comp. \cite{7}, (2.1) -- (2.7),
\item[(D)] metamorphoses of the $\zeta$-oscillating systems>
\begin{itemize}
\item[(a)] first, the notion of metamorphosis of an $\zeta$-oscillating multiform, \cite{4},
\item[(b)] after, the notion of metamorphosis of a quotient of two $\zeta$-oscillating multiforms, \cite{5}, \cite{6},
\end{itemize}
\item[(E)] $\mcal{Z}_{\zeta,Q^2}$-transformation (or device), \cite{7}.
\end{itemize}
Next, we have introduced (see \cite{8}) the following notions:
\begin{itemize}
\item[(F)] functional depending $\zeta$-oscillating systems, \cite{8}, beginning of the section 2.4,
\item[(G)] interaction between corresponding $\zeta$-oscillating systems, \cite{8}, Definition 4.
\end{itemize}

\subsection{}

In this paper we begin with the following set of functions
\be \label{1.1}
\begin{split}
& f_m(t)\in\tilde{C}_0[T,T+U],\ U=o\left(\frac{T}{\ln T}\right),\ T\to\infty, \\
& m=1,\dots,M,\ M\in \mbb{N},
\end{split}
\ee
$M$ being arbitrary and then fixed, where
\bdis
f_m\in \tilde{C}_0[T,T+U] \ \Leftrightarrow \ f_m\in C[T,T+U] \ \wedge \ f_m\not=0.
\edis
Next, by application of the operator $\hat{H}$ (introduced in \cite{8}, (3.6)) on elements (\ref{1.1}) we obtain the following vector-valued functions
\be \label{1.2}
\begin{split}
& \hat{H}f_m=(\alpha_0^{m,k_m},\alpha_1^{m,k_m},\dots,\alpha_{k_m}^{m,k_m},\beta_1^{k_m},\dots,\beta_{k_m}^{k_m}), \\
& m=1,\dots,M,\ 1\leq k_m\leq k_0,\ k_0\in\mbb{N},
\end{split}
\ee
where $k_0$ is arbitrary and fixed. Simultaneously, we obtain by our algorithm (short survey on this can be found in \cite{8}, (3.1) -- (3.11)) also the following set of
$\zeta$-factorization formulas
\be \label{1.3}
\begin{split}
& \prod_{r=1}^{k_m}\left|\frac{\zeta\left(\frac 12+i\alpha_r^{m,k_m}\right)}{\zeta\left(\frac 12+i\beta_r^{k_m}\right)}\right|^2\sim
E_m(U,T)F_m[f_m(\alpha_0^{m,k_m})],\\
& L\to\infty,\ m=1,\dots,M.
\end{split}
\ee

Now, we will suppose that after the finite number of stages of crossbreeding (every member of (\ref{1.3}) participates in this) in the set (\ref{1.3}) - that is: after the finite number of
elimination of the external functions
\be \label{1.4}
E_m(U,T),\ m=1,\dots,M
\ee
from the set (\ref{1.3}) - we obtain the following complete hybrid formula (i.e. the result of complete elimination of elements of (\ref{1.4})):
\be \label{1.5}
\begin{split}
& \mcal{F}\left\{\prod_{r=1}^{k_1}(\dots),\prod_{r=1}^{k_2}(\dots),\dots,\prod_{r=1}^{k_M}(\dots),F_1[f_1(\alpha_0^{1,k_1})],F_2[f_2(\alpha_0^{2,k_2})],\dots,F_M[f_M(\alpha_0^{M,k_M})]\right\} \\
& =1+o\left(\frac{\ln\ln T}{\ln T}\right)\sim 1,\ T\to\infty.
\end{split}
\ee

\begin{remark}
Here, we may put, of course,
\bdis
[T,T+U] \rightarrow [L,L+U], [\pi L,\pi L+U],\dots,\ L\in\mbb{N}.
\edis
\end{remark}

\begin{remark}
Let
\bdis
\mcal{F}_1\sim 1,\ T\to\infty
\edis
be the complete hybrid formula for the second set $\overline{(1.1)}$ on the segment $[T,T+U]$, where
\be \label{1.6}
\overline{(1.1)}\not=  (1.1).
\ee
Then (see (\ref{1.5}))
\be \label{1.7}
\mcal{F}\mcal{F}_1\sim 1,\ T\to\infty.
\ee
But, of course, (\ref{1.7}) is not the complete hybrid formula for the set
\bdis
(1.1) \cup \overline{(1.1)}
\edis
since lack of non-empty set of crossbreeding between sets (1.3) and (1.3).
\end{remark}

Now, we see that complete hybrid formula (\ref{1.5}) (it is simultaneously the interaction formula, comp. (G)) expresses the functional dependence (comp. (F)) of the set of vector-valued functions
(\ref{1.2}). Otherwise, the formula (\ref{1.5}) is a constraint on the set (\ref{1.2}), (comp. \cite{10}, Remark 9).

Consequently, the back-projection of this functional dependence of the set (\ref{1.2}) into the generating set (\ref{1.1}) leaves us to the following

\begin{mydef2}
We will call the subset
\bdis
\{ f_1(t),f_2(t),\dots,f_M(t)\},\ t\in  [T,T+U]
\edis
of the real continuous functions (comp. (\ref{1.1})) for which there is the complete hybrid formula (\ref{1.5}) as the family of $\zeta$-kindred functions.
\end{mydef2}

\subsection{}

In this paper, we prove that the following subset of the set of real continuous functions:
\be \label{1.8}
\begin{split}
& \{ \cos^2t,\sin^2t \},\ t\in [\pi L,\pi L+U],\ U\in (0,\pi/2),\ L\to\infty,
\end{split}
\ee
next,
\be \label{1.9}
\begin{split}
& \{ (t-L)^\Delta,(t-L)^{\Delta_1},\dots,(t-L)^{\Delta_n} \},\ t\in [ L, L+U],\ U\in (0,a),\ a<1,
\end{split}
\ee
where
\bdis
\Delta=\Delta_1+\dots +\Delta_n, \ \Delta>\Delta_1>\dots>\Delta_n>0,
\edis
and more
\be \label{1.10}
\begin{split}
& \left\{ \sum_{l=1}^n (t-L)^{\Delta_l},(t-L)^{\Delta_1},\dots,(t-L)^{\Delta_n}\right\}, \\
& t\in [L,L+U],\ U\in (0,a],\ a<1,\ L\to\infty, \\
& \Delta_l>0,\ \Delta_l\not=\Delta_k,\ l\not=k,\ 1\leq k,l\leq n
\end{split}
\ee
are the first families of $\zeta$-kindred functions.

\begin{remark}
The selection of families (\ref{1.8}) -- (\ref{1.10}) represents the completely new type of results in the theory of Riemann's zeta-function and, simultaneously, also in the
theory of real continuous functions.
\end{remark}

Finally, we notice explicitly, that also the results of this paper constitute the generic complement to the Riemann's functional equation on the critical line (comp. \cite{6}).

\section{The first family of $\zeta$-kindred trigonometric functions}

\subsection{}

If we put in \cite{8}, (4.1) -- (4.10)
\bdis
\mu=0,\ \tilde{C}\rightarrow \tilde{C}_0
\edis
(comp. \cite{9}, Definition) then we obtain the following lemmas.

\begin{mydef1}
For the function
\be \label{2.1}
f_1(t)=\sin^2t\in\tilde{C}_0[\pi L,\pi L+U],\ U\in (0,\pi/2)
\ee
there are vector-valued functions
\be \label{2.2}
(\alpha_0^{1,k_1},\alpha_1^{1,k_1},\dots,\alpha_{k_1}^{1,k_1},\beta_1^{k_1},\dots,\beta_{k_1}^{k_1})
\ee
such that the following factorization formula
\be \label{2.3}
\begin{split}
& \prod_{r=1}^{k_1}\left|\frac{\zeta\left(\frac 12+i\alpha_r^{1,k_1}\right)}{\zeta\left(\frac 12+i\beta_r^{k_1}\right)}\right|^2\sim \\
& \left(\frac 12-\frac 12\frac{\sin U}{U}\cos U\right)\frac{1}{\sin^2(\alpha_0^{1,k_1})},\ L\to\infty
\end{split}
\ee
holds true, where
\be \label{2.4}
\begin{split}
& \alpha_r^{1,k_1}=\alpha_r(U,L,k_1;f_1),\ r=0,1,\dots,k_1, \\
& \beta_r^{k_1}=\beta_r(U,L,k_1),\ r=1,\dots,k_1, \\
& 0<\alpha_0^{1,k_1}-\pi L<U.
\end{split}
\ee
\end{mydef1} 

\begin{mydef52} 
For the function
\be \label{2.5}
f_2(t)=\cos^2t\in\tilde{C}_0[\pi L,\pi L+U],\ U\in (0,\pi/2)
\ee
there are vector-valued functions
\be \label{2.6}
(\alpha_0^{2,k_2},\alpha_1^{2,k_2},\dots,\alpha_{k_2}^{2,k_2},\beta_1^{k_2},\dots,\beta_{k_2}^{k_2})
\ee
such that the following $\zeta$-factorization formula
\be \label{2.7}
\begin{split}
	& \prod_{r=1}^{k_2}\left|\frac{\zeta\left(\frac 12+i\alpha_r^{2,k_2}\right)}{\zeta\left(\frac 12+i\beta_r^{k_2}\right)}\right|^2\sim \\
	& \left(\frac 12-\frac 12\frac{\sin U}{U}\cos U\right)\frac{1}{\cos^2(\alpha_0^{2,k_2})},\ L\to\infty
\end{split}
\ee
holds true, where
\be \label{2.8}
\begin{split}
	& \alpha_r^{2,k_2}=\alpha_r(U,L,k_2;f_2),\ r=0,1,\dots,k_2, \\
	& \beta_r^{k_1}=\beta_r(U,L,k_2),\ r=1,\dots,k_2, \\
	& 0<\alpha_0^{2,k_2}-\pi L<U.
\end{split}
\ee
\end{mydef52} 

\subsection{}  

Crossbreeding between the $\zeta$-factorization formulae (\ref{2.3}) and (\ref{2.7}) -- one stage is sufficient in this case -- gives the following 

\begin{mydefCHF1} 
\be \label{2.9} 
\begin{split} 
&  \cos^2(\alpha_0^{2,k_2})\prod_{r=1}^{k_2}\left|\frac{\zeta\left(\frac 12+i\alpha_r^{2,k_2}\right)}{\zeta\left(\frac 12+i\beta_r^{k_2}\right)}\right|^2+ \\ 
& + \sin^2(\alpha_0^{1,k_1}) \prod_{r=1}^{k_1}\left|\frac{\zeta\left(\frac 12+i\alpha_r^{1,k_1}\right)}{\zeta\left(\frac 12+i\beta_r^{k_1}\right)}\right|^2\sim 1,\ L\to\infty . 
\end{split} 
\ee 
\end{mydefCHF1}  

Now, we obtain from (\ref{2.9}) by Definition the following 

\begin{mydef11} 
	The subset 
\be \label{2.10} 
\begin{split} 
& \{ \cos^2t,\sin^2t\},\ t\in [\pi L,\pi L+U],\ U\in (0,\pi/2),\ L\to\infty 
\end{split} 
\ee     
is the family of $\zeta$-kindred elements in the class of trigonometric functions. 
\end{mydef11}  

\begin{remark} 
	The formula (\ref{2.9}) has already been obtained in our paper \cite{8}, (2.6). However, in the present paper the same one is playing the role of complete hybrid 
	formula. Of course, the set 
	\bdis 
	\{ \cos^2t,\sin^2t\};\ \cos^2t+\sin^2t=1 
	\edis  
	is the known family of \emph{school-kindred} functions. 
\end{remark} 

\section{The first family of $\zeta$-kindred real power functions} 

\subsection{} 

The following lemma holds true (see \cite{9}, (2.4), (2.5)). 
\begin{mydef53} 
For the function 
\be \label{3.1} 
\bar{f}_\Delta(t,L)=\bar{f}_\Delta(t)=(t-L)^\Delta\in \tilde{C}_0[L,L+U],\ U\in (0,a],\ a<1,\ \Delta>0 
\ee  
there are vector-valued  functions 
\be \label{3.2} 
(\bar{\alpha}_0^{\Delta,\bar{k}_\Delta},\bar{\alpha}_1^{\Delta,\bar{k}_\Delta},\dots,\bar{\alpha}_{\bar{k}_\Delta}^{\Delta,\bar{k}_\Delta},
\bar{\beta}_1^{\bar{k}_\Delta},\dots,\bar{\beta}_{\bar{k}_\Delta}^{\bar{k}_\Delta}) 
\ee 
such that the following $\zeta$-factorization formula 
\be \label{3.3}
\begin{split}
	& \prod_{r=1}^{\bar{k}_\Delta}\left|\frac{\zeta\left(\frac 12+i\bar{\alpha}_r^{\Delta,\bar{k}_\Delta}\right)}{\zeta\left(\frac 12+i\bar{\beta}_r^{\bar{k}_\Delta}\right)}\right|^2\sim 
	\frac{1}{1+\Delta}\left(\frac{U}{\bar{\alpha}_0^{\Delta,\bar{k}_\Delta}-L}\right)^\Delta,\ L\to\infty
\end{split}
\ee 
holds true, where 
\be \label{3.4}
\begin{split}
	& \bar{\alpha}_r^{\Delta,\bar{k}_\Delta}=\alpha_r(U,L,\bar{k}_\Delta;\bar{f}_\Delta),\ r=0,1,\dots,\bar{k}_\Delta, \\
	& \bar{\beta}_r^{\bar{k}_\Delta}=\beta_r(U,L,\bar{k}_\Delta),\ r=1,\dots,\bar{k}_\Delta, \\
	& 0<\bar{\alpha}_0^{\Delta,\bar{k}_\Delta}-\pi L<U.
\end{split}
\ee
\end{mydef53}  

\subsection{} 

Let 
\be \label{3.5} 
\Delta=\sum_{l=1}^n \Delta_l,\ \Delta>\Delta_1>\dots>\Delta_n>0,\ n\in\mbb{N} 
\ee  
(for every fixed $n$). Then we have by Lemma 3 that 
\be \label{3.6} 
\begin{split} 
& (1+\Delta_l)(\bar{\alpha}_0^{\Delta_l,k_l}-L)^{\Delta_l}\prod_{r=1}^{k_l}\left|\frac{\zeta\left(\frac 12+i\bar{\alpha}_r^{\Delta_l,k_l}\right)}{\zeta\left(\frac 12+i\bar{\beta}_r^{k_l}\right)}\right|^2
\sim U^{\Delta_l}, \\ 
& L\to\infty,\ \bar{k}(\Delta_l)\rightarrow k_l,\ \bar{\alpha}_r^{\Delta,\bar{k}(\Delta_l)}\rightarrow \bar{\alpha}_r^{\Delta_l,k_l},\ \dots 
\end{split} 
\ee  
i.e. we have the $n$-analogues of (\ref{3.4}). 

Now, after the $n$-stages of crossbreeding between the $\zeta$-factorization formulae (\ref{3.3}) and (\ref{3.6}) where, of course, 
\be 
U^{\Delta_1+\dots+\Delta_n}=U^{\Delta_1}\cdot \dots\cdot U^{\Delta_n}, 
\ee 
we obtain the following result.  

\begin{mydefCHF2}
\be \label{3.8} 
\begin{split} 
&  \prod_{r=1}^{\bar{k}_\Delta}\left|\frac{\zeta\left(\frac 12+i\bar{\alpha}_r^{\Delta,\bar{k}_\Delta}\right)}{\zeta\left(\frac 12+i\bar{\beta}_r^{\bar{k}_\Delta}\right)}\right|^2 \sim \\ 
& \sim \frac{1}{1+\Delta}\prod_{l=1}^n (1+\Delta_l)\frac{1}{(\bar{\alpha}_0^{\Delta,\bar{k}_\Delta}-L)^\Delta}
\prod_{l=1}^n(\bar{\alpha}_0^{\Delta_l,k_l}-L)^{\Delta_l}\times \\ 
& \times \prod_{l=1}^n\prod_{r=1}^{k_l}\left|\frac{\zeta\left(\frac 12+i\bar{\alpha}_r^{\Delta_l,k_l}\right)}{\zeta\left(\frac 12+i\bar{\beta}_r^{k_l}\right)}\right|^2,\ L\to\infty, \\ 
& 1\leq \bar{k}_\Delta,k_1,\dots,k_n\leq k_0;\ k_1=\bar{k}(\Delta_1),\dots,
\end{split} 
\ee  
here, of course, 
\bdis 
(\dots)\sim [\dots] \ \Leftrightarrow \ 1\sim \frac{(\dots)}{[\dots ]};\ (\dots),[\dots]\not=0. 
\edis 
\end{mydefCHF2} 

\begin{remark} 
	The symbol (\ref{3.8}) contains the set of 
	\bdis 
	(k_0)^{n+1} 
	\edis  
	formulas. For example, in the case $k_0=100$ and $n=99$ this number is equal to 
	\bdis 
	10^{200} . 
	\edis 
\end{remark}   

Consequently, from (\ref{3.8}) we obtain by Definition the following   

\begin{mydef12}
The subset 
\be \label{3.9} 
\begin{split} 
&	\left\{ (t-L)^\Delta, (t-L)^{\Delta_1},\dots,(t-L)^{\Delta_n}\right\}, \\ 
& t\in [L,L+U], U\in (0,a],\ a<1,\ L\to\infty, \\ 
& \Delta=\sum_{l=1}^n \Delta_l,\ \Delta>\Delta_1>\dots>\Delta_n>0 
\end{split}
\ee  
is the family of $\zeta$-kindred elements in the class of real power functions. 
\end{mydef12}  

\begin{remark}
	The formula (\ref{3.8}) has already been obtained in our paper \cite{9}, (3.3). However, in the present paper this one is playing the role of the complete hybrid formula.
\end{remark} 

\section{The second family of $\zeta$-kindred real power functions} 

\subsection{} 

The following lemma holds true (see \cite{9}, (4.2)). 

\begin{mydef54}
	For the function 
	\be \label{4.1} 
	\begin{split}
		& \tilde{f}(t)=\tilde{f}(t;\Delta_1,\dots,\Delta_n,L)=\sum_{l=1}^n (t-L)^{\Delta_l}\in \tilde{C}_0[L,L+U],\\ 
		& U\in (0,a],\ a<1,\ \Delta_l>0, \\ 
		& \Delta_l\not=\Delta_k,\ l\not= k;\ l,k=1,\dots,n 
	\end{split} 
    \ee  
    there are vector-valued functions 
    \be \label{4.2} 
    (\tilde{\alpha}_0,\tilde{\alpha}_1,\dots,\tilde{\alpha}_k,\tilde{\beta}_1,\dots,\tilde{\beta}_k),\ 1\leq k\leq k_0 
    \ee  
    such that the following $\zeta$-factorization formula 
    \be \label{4.3} 
    \prod_{r=1}^k
    \left|\frac{\zeta\left(\frac 12+i\tilde{\alpha}_r\right)}{\zeta\left(\frac 12+i\tilde{\beta}_r\right)}\right|^2\sim 
    \frac
    {\sum_{l=1}^n\frac{1}{1+\Delta_l}U^{\Delta_l}}
    {\sum_{l=1}^n (\tilde{\alpha}_0-L)^{\Delta_l}},\ L\to\infty 
    \ee  
    holds true, where 
    \be \label{4.4} 
    \begin{split}
    	& \tilde{\alpha}_r=\alpha_r(U,L,\Delta_1,\dots,\Delta_n,k),\ r=0,1,\dots,k, \\ 
    	& \tilde{\beta}_r=\beta_r(U,L,k),\ r=1,\dots,k. 
    \end{split} 
    \ee 
\end{mydef54} 

\subsection{} 

Next, we have (comp. (\ref{3.6})) the following formulas 
\be \label{4.5} 
\begin{split}
& (\alpha_0^{\Delta_l,k_l}-L)^{\Delta_l}\prod_{r=1}^k 
\left|\frac{\zeta\left(\frac 12+i\alpha_r^{\Delta_l,k_l}\right)}{\zeta\left(\frac 12+i\beta_r^{k_l}\right)}\right|^2\sim \frac{1}{1+\Delta_l}U^{\Delta_l}, \\ 
& \Delta_l>0,\ 1\leq k_l\leq k_0,\ l=1,\dots,n. 
\end{split} 
\ee  
Now we obtain after $n$-stages of crossbreeding between the $\zeta$-factorization formulas (\ref{4.3}) and (\ref{4.5}) the following result (see \cite{9}, (4.4)).  

\begin{mydefCHF3}
\be \label{4.6} 
\begin{split}
	& \prod_{r=1}^k\left|\frac{\zeta\left(\frac 12+i\tilde{\alpha}_r\right)}{\zeta\left(\frac 12+i\tilde{\beta}_r\right)}\right|^2\sim \\ 
	& \sim \frac{1}{\sum_{l=1}^n (\tilde{\alpha}_0-L)^{\Delta_l}} \sum_{l=1}^n (\alpha_0^{\Delta_l,k_l}-L)^{\Delta_l}
	\prod_{r=1}^n\left|\frac{\zeta\left(\frac 12+i\alpha_r^{\Delta_l,k_l}\right)}{\zeta\left(\frac 12+i\beta_r^{k_l}\right)}\right|^2. 
\end{split}
\ee 
\end{mydefCHF3} 

Consequently, we obtain from (\ref{4.6}) by Definition the following 

\begin{mydef13}
The subset 
\be \label{4.7} 
\begin{split}
	& \left\{ \sum_{l=1}^n(t-L)^{\Delta_l},\ (t-L)^{\Delta_1},\dots,(t-L)^{\Delta_l}\right\}, \\ 
	& t\in [L,L+U],\ U\in (0,a],\ a<1,\ L\to\infty , \\ 
	& \Delta_l>0,\ \Delta_l\not=\Delta_k,\ l\not=k,\ 1\leq k,l\leq n 
\end{split}
\ee  
is the family of $\zeta$-kindred elements in the class of real power functions. 
\end{mydef13} 

\section{Concluding remarks} 

\subsection{} 

The first remarks are connected with the operator $\hat{H}$ that we have defined in the paper \cite{8}, (see Definition 2 and Definition 5; comp. also \cite{9}, Definition). In this 
direction we have used the following notation 
\be \label{5.1} 
\begin{split}
&	\forall\- f(t)\in \tilde{C}_0[T,T+U] \rightarrow  \hat{H}f(t)=(\alpha_0,\alpha_1,\dots,\alpha_k,\beta_1,\dots,\beta_k),\ k=1,\dots,k_0 
\end{split}
\ee 
for every fixed $k$, (see \cite{10}, section 5). 

\begin{remark}
	Here we have to use more exact notation 
\be \label{5.2} 
(\alpha_0,\alpha_1,\dots,\alpha_k,\beta_1,\dots,\beta_k) \rightarrow (\alpha_0^k,\alpha_1^k,\dots,\alpha_k^k,\beta_1^k,\dots,\beta_k^k)
\ee 
Namely, we have defined in \cite{10}, section 5, the $\hat{H}$ operator as matrix-valued operator with the following inexact notation based on (\ref{5.1}): 
\be \label{5.3} 
f(t)\xrightarrow{\hat{H}} 
\begin{pmatrix}
	\alpha_0 & \alpha_1 & \beta_1 & 0 & \dots\\ 
	\alpha_0 & \alpha_1 & \alpha_2 & \beta_1 & \beta_2 & 0 & \dots \\ 
	\vdots \\ 
	\alpha_0 & \alpha_1 & \dots & \alpha_{k_0} & \beta_1 & \beta_2 & \dots & \beta_{k_0}
\end{pmatrix}_{k_0\times (2k_0+1)} . 
\ee  
However, the exact notation (instead of (\ref{5.3})) is the following one (comp. (\ref{5.2}))  
\be \label{5.4} 
f(t)\xrightarrow{\hat{H}} 
\begin{pmatrix}
	\alpha_0^1 & \alpha_1^1 & \beta_1^1 & 0 & \dots\\ 
	\alpha_0^2 & \alpha_1^2 & \alpha_2^2 & \beta_1^2 & \beta_2^2 & 0 & \dots \\ 
	\vdots \\ 
	\alpha_0^{k_0} & \alpha_1^{k_0} & \dots & \alpha_{k_0}^{k_0} & \beta_1^{k_0} & \beta_2^{k_0} & \dots & \beta_{k_0}^{k_0}
\end{pmatrix}_{k_0\times (2k_0+1)} . 
\ee  
\end{remark} 

\begin{remark}
	Now, it is clear that the inexact notation (\ref{5.3}) may suggest erroneous impression that the elements of the column 
	\bdis 
	\begin{pmatrix} \alpha_0 \\ \alpha_0 \\ \vdots \\ \alpha_0 \end{pmatrix} , 
	\edis 
	for example, are necessarily mutually equal. 
\end{remark}

\subsection{} 

The second remark is connected with the asymptotic form: 
\begin{itemize}
	\item[(a)] of the $\zeta$-factorization formulas (\ref{1.3}), (\ref{2.3}), (\ref{2.7}), (\ref{3.3}), (\ref{4.1}),  
	\item[(b)] of the complete hybrid formulas (\ref{1.5}), (\ref{2.9}), (\ref{3.8}), (\ref{4.6}). 
\end{itemize} 

\begin{remark}
	The following is true: if we use from beginning the exact factorization formula 
	\bdis 
	\prod_{r=1}^k \frac{\tilde{Z}^2(\alpha_r)}{\tilde{Z}^2(\beta_r)}=\frac{H(T,U;f)}{f(\alpha_0)},\ T\to\infty 
	\edis  
	instead of the asymptotic $\zeta$-factorization formula 
	\bdis 
	\prod_{r=1}^k\left|\frac{\zeta\left(\frac 12+i\alpha_r\right)}{\zeta\left(\frac 12+i\beta_r\right)}\right|^2= 
	\left\{ 1+\mcal{O}\left(\frac{\ln\ln T}{\ln T}\right)\right\} \frac{H(T,U;f)}{f(\alpha_0)}\sim \frac{H(T,U;f)}{f(\alpha_0)},\ T\to\infty 
	\edis  
	(see short survey of our algorithm for generating $\zeta$-factorization formulas in \cite{8}, (3.7), (3.8)), then we obtain the exact $\zeta$-factorization formulas and also the exact 
	complete hybrid formula instead of these in (a) and (b). 
\end{remark}

\end{document}